\documentclass[11pt,leqno]{article}
\usepackage{amsmath,amsthm}
\usepackage {latexsym}
\usepackage{amssymb}
\newcommand{\les}{\lesssim}
\newcommand{\bea}{\begin{eqnarray}}\newcommand{\eea}{\end{eqnarray}}
\newcommand{\beq}{\begin{equation}}\newcommand{\ee}{\end{equation}}
\def\pa{\partial}

\newtheorem{theorem}{Theorem}[section]\newtheorem{lemma}[theorem]{Lemma}

\theoremstyle{remark}\newtheorem{remark}[theorem]{Remark}

\def\bm{\left( \begin{array}{cc}}\def\endm{\end{array}\right)}
\newcommand{\eq}{\end{equation}}
\def\pa{\partial}
\def \rectangle#1#2{\hbox{\vrule\vbox to #2
              {\hrule\hbox to #1{\hfil}\vfil\hrule}\vrule}}

\def\ro{{\varrho_{\,}}}

\numberwithin{equation}{section}

\begin{document}
\author {Hans Lindblad\thanks{Part of this work was done while H.L. was a Member of the
Institute for  Advanced Study, Princeton, supported by the NSF
grant DMS-0111298 to the Institute. H.L. was also partially
supported by the NSF Grant DMS-0200226. } \, and Avy
Soffer\thanks{Also a member of the Institute of Advanced
Study, Princeton.Supported in part by NSF grant DMS-0100490.}\\
University of California at San Diego and Rutgers University}
\title {A remark on asymptotic completeness for the critical nonlinear Klein-Gordon equation} \maketitle
\begin{abstract}
We give a short proof of asymptotic completeness and global
existence for the cubic Nonlinear Klein-Gordon equation in one
dimension. Our approach to dealing with the long range behavior of
the asymptotic solution is by reducing it, in hyperbolic
coordinates to the study of an ODE. Similar arguments extend to
higher dimensions and other long range type nonlinear problems.
\end{abstract}
\section{Introduction} We are interested in the completeness problem for the critical nonlinear Klein-Gordon in one space dimension:
 \beq \label{eq:nlkg}
 \square v+ v=-\beta v^3,\qquad \beta\geq 0
 \eq
 where $\square=\pa_t^2-\pa_x^2$.

 Such problem appear naturally in the study of some nonlinear
 dynamical problems of mathematical physics, among them radiation theory \cite{F}, general relativity \cite{L-R}
 , the
 scattering and stability of kinks, vortices and other coherent
 structures. The long range nature of the scattering, coupled with
 the nonlinearity poses a challenge to scattering theory.
 In the Linear case, one can compute in advance the asymptotic
 corrections to the solution due to the long range nature of the
 solution. That is not possible in the nonlinear case. It may also break global existence \cite{H2}.
 Since the asymptotic behavior is not free other tools are needed
 to prove scattering and global existence, e.g. normal forms \cite{Sim} for dealing with $u^2$ terms
 \cite{Sh,D1}. Inspired by Delort's
 ground breaking work
 \cite{D1}, here we see
 that one can nevertheless get the correct asymptotic behavior
 from that of an ODE.
 Recall first that a solution of
the linear Klein-Gordon, i.e. $\beta=0$, is asymptotically given by
(as $\ro$  tends to infinity)
 \beq \label{eq:linearasymptotic}
 u(t,x)\sim \ro^{-1/2} e^{i\ro}\, a(x/\!\ro)+
 \ro^{-1/2} e^{-i\ro}\, \overline{a(x/\!\ro)},
 \qquad \text{where}\quad \ro=(t^2-|x|^2)^{1/2}\geq 0.
 \eq
Here $a(x/\!\ro)=(t/\!\ro)
\widehat{u}_+(-x/\!\ro)=\sqrt{1+x^2/\ro^2}\,
\widehat{u}_+(-x/\!\ro)$, where $\widehat{u}_+(\xi)=\int u_+(x) \,
e^{-ix \xi}\, dx$ denotes the Fourier transform with respect to $x$
only, $\widehat{u}_+=(\widehat{u}_0-i(|\xi|^2+1)^{-1/2}
\widehat{u}_1)/2$, where $u_0=u\big|_{t=0}$ and $u_1=\pa_t
u\big|_{t=0}$. Here the right hand side is to be interpreted as $0$
outside the light cone, when $|x|>t$. \eqref{eq:linearasymptotic}
can be proven using stationary phase, see e.g. \cite{H1}, where a
complete asymptotic expansion into negative powers of $\ro$ was
given. One can also get precise estimates by vector-field
multipliers \cite{K1}. Recently, Delort\cite{D1} showed that
\eqref{eq:nlkg} with small initial data have a global solution with
asymptotics of the form
 \beq
\label{eq:nonlinearasymptotic}
 v(t,x)\sim \ro^{-1/2} e^{i\phi_0(\ro\!,\,x/\ro)}\, a(x/\!\ro)+
 \ro^{-1/2} e^{-i\phi_0(\ro\!,\, x/\ro)}\, \overline{a(x/\!\ro)},
 \qquad
 \phi_0(\ro,x/\!\ro)=\ro+\frac{3}{8}\beta |a(x/\!\ro)|^2\,\ln{\ro}
 \eq
 Delort's proof is for more general equations but it is rather involved ;
 furthermore there is an application of Gronwall's lemma,
  used to obtain the estimate (3.27) from (3.26)
  that seem to require an additional argument. This is because the
  integrand in (3.26) depends explicitly on t.
 In \cite{L-S1} we considered the inverse problem of scattering, we showed that
 for any given asymptotic expansion of the above form (1.3) there is a
 solution agreeing with it at infinity.
  The purpose of this note is to present a simple proof of the completeness for
 the special case considered here.

For technical reasons we will give initial data of compact support
$|x|\leq 1$ in the support of $u_0$ and $u_1$ when $t=2$:
$$
v(2,x)=\varepsilon u_0(x), \qquad \pa_t v(2,x)=\varepsilon u_1(x).
$$
so that solution is supported inside the forward light cone
$t-|x|\geq 1$, when $t\geq 2$.
\begin{remark}
The above restriction on the initial data makes it difficult to
give an explicit description of the space in which the asymptotic
solutions live. An extension of the proof to more general data,
such as the Schwartz class is therefore desirable, and probably
can be achieved following arguments similar to\cite{S-Taf}
\end{remark}

 For the proof
we start by introducing the hyperbolic coordinates
$$
\ro^2=t^2-x^2, \quad t=\ro\cosh y,\quad x=\ro\sinh y,
$$
or \beq\label{eq:yweight}
 e^{2|y|}=\frac{t+|x|}{t-|x|},\qquad
\ro^2=t^2-x^2 \eq Then
$$
\square+1=\pa_\ro^2-\ro^{-2}\pa_y^2+\ro^{-1}\pa_\ro+1
$$
and with
$$
v(t,x)=\ro^{-1/2}V(\ro,y)
$$
we get
$$
(\square+1)v(t,x)=\ro^{-1/2}\Big(
\pa_\ro^2+1-\ro^{-2}(\pa_y^2-\tfrac{1}{4})\Big) V(\ro,y).
$$

Hence in these coordinates \eqref{eq:nlkg} becomes the following
equation for $V=\ro^{1/2}v$:
 \beq \label{eq:PDE} \Psi(V)\equiv\pa_\ro^2
V+\Big(1+\frac{\beta}{\ro}V^2+\frac{1}{4\ro^2}\Big)V-\frac{1}{\ro^{2}}\pa_y^2
V=0 \eq
 We are therefore led to first studying the ODE
 \beq\label{eq:ODE}
L(g)\equiv \ddot{g}+\Big(1+\frac{\beta}{\ro}
g^2+\frac{1}{4\ro^2}\Big)g=F \eq

The proof as usual consists of a decay estimate ($L^\infty$
estimate) of a lower number of derivatives and energy estimate
($L^2$ estimate) of a higher number of derivatives. The $L^2$
estimate uses the decay estimate $C\varepsilon \rho^{-1/2}$ and
allows for growing energies, like $Ct^\varepsilon$. Since the
energies are growing we don't directly get back the $L^\infty$
estimate from the energy estimates but instead we get some weaker
decay estimates using some weighted Sobolev lemmas. Using these
weaker decay estimates we can get further decay from the asymptotic
equation, i.e. the ODE above, using that the term $\rho^{-2} \pa_y^2
V$ has a negative power of $\rho$ so there is some room and we can
use the weak decay estimates to estimate this term, c.f. \cite{L1},
\cite{L-R} The energies we will use will just be the energies on
hyperboloids in the coordinates\eqref{eq:yweight}

\section{The sharp decay estimate}

\begin{lemma}\label{decayode} Suppose that
\beq\label{eq:ODE2} \ddot{g}+\Big(1+\frac{\alpha}{\ro^{1/2}}
g+\frac{\beta}{\ro} g^2+\frac{1}{4\ro^2}\Big)g=F \eq Then

\beq
 |\dot{g}(\ro)|+|g(\ro)|\leq 2\Big( |\dot{g}(1)|+|g(1)|+
|g(1)|^2 +\int_1^\ro |F(\tau)|\, d\tau\Big)
\eq
 if $\alpha=0$ and
$\beta\geq 0$, or the right hand side is sufficiently small.
\end{lemma}
\begin{proof}
 Multiplying \eqref{eq:ODE2} by  $2\dot{g}$
 we see that
 $$\frac{d}{d\ro}
 \Big( \dot{g}^2+g^2+\frac{2\alpha}{3\ro^{1/2}} g^3 +\frac{\beta}{2\ro} g^4 + \frac{1}{4\ro^2} g^2 \Big)=2F\dot{g} -\frac{\alpha}{3\ro^{3/2}} g^3 -\frac{\beta}{2\ro^2} g^4 - \frac{1}{2\ro^3} g^2
 $$
 It therefore follows that if $\alpha=0$ and $\beta\geq 0$ then
$|2M\frac{dM}{d\rho}|\leq |2F\dot{g}|$, so
\beq
\Big|\frac{d}{d\ro} M(\ro)\Big|\leq F(\ro),\qquad\text{with}\quad
M=\Big( \dot{g}^2+g^2+\frac{2\alpha}{3\ro^{1/2}} g^3+\frac{\beta}{2\ro} g^4 + \frac{1}{4\ro^2}g^2
\Big)^{1/2}
\eq
\end{proof}
Applying the above lemma to the equation (1.5)\beq \label{eq:PDE2}
\pa_\ro^2
V+\Big(1+\frac{\beta}{\ro}V^2+\frac{1}{4\ro^2}\Big)V=\frac{1}{\ro^{2}}\pa_y^2
V \eq and noting that $|\dot{V}|+|V|\leq C_0\varepsilon $ when
$\ro=1$ by the support assumptions on initial data, we conclude
that \beq |\pa_\ro V(\rho,y)|+|V(\ro,y)|\leq C\varepsilon
+\int_1^\ro \tau^{-2} |\pa_y^2 V(\tau,y)|\, d\tau \eq Which gives
the bound \beq |\pa_\ro V(\rho,y)|+|V(\ro,y)|\leq C_2\varepsilon
\eq if we can prove that for some $\delta>0$; \beq |\pa_y^2
V(\ro,y)|\leq C_1 \varepsilon (1+\ro)^{1-\delta} \eq

\section{Energies on hyperboloids and the weak decay estimate}
 By section 7.6 in H\"ormander \cite{H1}, if
 $$
 E(\ro)^2=\int_{H_\ro} u_t^2+u_x^2+2\frac{x}{t} u_t u_x+ u^2 \, dx
 $$
 where $H_\ro=\{(t,x);\, t^2-x^2=\ro^2\}$
 and $G_{\ro_1,\ro_2}=\{(t,x);\, \ro_2\geq t^2-x^2\geq \ro_1\}$ then
(recall $\frac{|x|}{t}\leq 1\, \frac{\rho}{t}\leq 1$)
 \beq
 \frac{d}{d\ro} E(\ro)^2=2\int_{H_\ro} F \, u_t\frac{\ro}{t} \, dx\leq
2\Big(\int_{H_\ro} |F|^2\, dx\Big)^{1/2} E(\ro)
 \eq
and hence

$$E(\ro)\leq E(1)+\int_1^\ro \Big(\int_{H_\tau} |F|^2\,
dx\Big)^{1/2}\, d\tau.$$
We have
$$u_t^2+u_x^2+2\frac{x}{t} u_t u_x = \Big(u_x+\frac{x}{t}
u_t\Big)^2+\frac{\ro^2}{t^2} u_t^2 =\Big(u_t+\frac{x}{t}
u_x\Big)^2+\frac{\ro^2}{t^2} u_x^2.
 $$
Since $u_\ro=(tu_t+x u_x)/\ro$, $u_y=xu_t+tu_x$
and $t=\rho\cosh{y}$ we see that
\beq
u_t^2+u_x^2+2\frac{x}{t} u_t u_x
=\frac{1}{2\cosh^2{y}}\Big( u_\ro^2+\frac{u_y^2}{\ro^2}\Big)
+\frac{\ro^2}{2
t^2}\big(u_t^2+u_x^2\big)
\eq
 In order to change variables in the integral we think of $x=\ro\sinh y$ as a
 function of $y$ for $\ro$ fixed. Then $\pa x/\pa
 y\big|_{\ro=const}=\ro \cosh y=t \sim \ro e^{|y|}$.
Hence with $V=\ro^{1/2} u$ we have
\beq
E(\ro)^2\geq c\int_{H_\ro} \big(V_y^2\ro^{-2}+V_\ro^2\big) e^{-|y|}+V^2 e^{|y|}\,dy,\qquad \ro\geq 1.
\eq
Here we also used that $V_\ro=\ro^{1/2}u_\ro+V/2\ro$ and
$V^2\ro^{-2}\cosh^{-2}{y}=V^2t^{-2}\leq V^2$ if $t\geq 1$.

Moreover we have proven using (3.2), (3.3):

\begin{lemma} Suppose that
\beq
\pa_t^2 w-\pa_x^2 w+w=F
\eq
Then
\beq
\Big(\int_{H_\ro} V_y(\ro,y)^2\ro^{-2} e^{-|y|}+V(\ro,y)^2 e^{|y|}\, dy\Big)^{1/2}
\leq
C E(1)+\int_1^\ro \tau^{1/2} \Big(\int_{H_\tau} |F(\tau,y)|^2 e^{|y|}\, dy \Big)^{1/2}\, d\tau
\eq
\end{lemma}

\begin{lemma}\label{Energy2} Suppose that
\beq
\pa_\ro^2 W+\big(1+\frac{1}{4\ro^2}\big) W-\frac{1}{\ro^{2}}\pa_y^2 W=F
\eq
Let
\beq
\|F(\ro,\cdot)\|_{L^p(H_\ro)}
=\Big(\int_{H_\ro} |F(\ro,y)|^p\, e^{|y|}\, dy\Big)^{1/p}
\eq
Then
\beq\label{eq:inhomenergy1}
 \|W(\ro,\cdot)\|_{L^2(H_\ro)}\leq CE(1)
 +C\int_1^{\ro} \|F(\sigma,\cdot)\|_{L^2(H_\sigma)}\, d\sigma
 \eq
\end{lemma}
In the applications
\beq
\pa_\ro^2 V^{(k)}+\big(1+\frac{1}{4\ro^2}\big) V^{(k)}-\frac{1}{\ro^{2}}\pa_y^2 V^{(k)}=F^{(k)}
\eq
where
\beq
F^{(k)}= \beta \ro^{-1} \pa_y^k V^3,\qquad V^{(k)}=\pa_y^k V
\eq
We claim that
\beq
\|F^{(k)}(\ro,\cdot)\|_{L^2(H_\ro)}\leq C \ro^{-1} \|V(\ro,\cdot)\|_{L^\infty(H_\ro)}^2 \|V(\ro,\cdot)\|_{L^{2,k}(H_\ro)} ,\quad k=0,1,2,3,
\eq
where
\beq
\|F(\ro,\cdot)\|_{L^{p,k}(H_\ro)}
=\Big(\int_{H_\ro} \sum_{m\leq k}|\pa_y^m F(\ro,y)|^p\, e^{|y|}\, dy\Big)^{1/p}
\eq
In fact
\begin{align}
|F^{(0)}|&\leq C \ro^{-1}  |V|^3\\
|F^{(1)}|& \leq C \ro^{-1}  |V|^2 |\pa_y V|\\
|F^{(2)}|&\leq C \ro^{-1}\big(|V|^2\, |\pa_y^2 V| +|V|\, |\pa_y V|^2\big)\\
|F^{(3)}|&\leq C \ro^{-1} \big(|V|^2|\pa_y^3 V| +|V|\, |\pa_y V|\, |\pa_y^2 V|+ |\pa_y V|^3\big)
\end{align}
For $k=0,1$ (3.13),(3.14) are obvious and for $k\geq 2$
(3.15),(3.16) follow from interpolation:
\begin{lemma}

For $k\geq 2\, j\leq k$
\beq
\|V\|_{L^{2k/\!j, j}(H_\ro)}^{k/j}\leq C\|V\|_{L^\infty(H_\ro)}^{k/j-1} \| V\|_{L^{2,k}(H_\ro)}
\eq
\end{lemma}
\begin{proof} A proof of the standard interpolation without weights
just uses H\"older's inequality, which holds with weights, and
integration by parts, which produces lowers order terms included in
the norms.
\end{proof}
For $k=2$ we have $(j=1)$ \beq \|V(t,\cdot)\|_{L^{4,1}}^2\leq C
\|V(t,\cdot)\|_{L^\infty} \| V(t,\cdot)\|_{L^{2,2}}, \eq Similarly
for $k=3$ we have $j=(1,2)$
\begin{align}
\|V(t,\cdot)\|_{L^{6,1}}^3 &\leq
C \|V(t,\cdot)\|_{L^\infty}^{1/2} \|V(t,\cdot)\|_{L^{2,3}},\\
\| V(t,\cdot)\|_{L^{3,2}}^{3/2} &\leq
C \|V(t,\cdot)\|_{L^\infty}^2 \| V(t,\cdot)\|_{L^{2,3}},
\end{align}
Assuming the bound \beq \|V(\ro,\cdot)\|_{L^\infty(H_\ro)}\leq
C_0\varepsilon \, \rho\geq 1 \eq with a constant independent of
$\ro$ we have hence proven that using (3.6)-(3.8)
\beq\label{eq:inhomenergy1}
 \|V(\ro,\cdot)\|_{L^{2,k}(H_\ro)}\leq CE(1)
 +\int_1^{\ro} C\varepsilon^2 \sigma^{-1}\|V(\sigma,\cdot)\|_{L^{2,k}(H_\sigma)}\, d\sigma,\quad k=0,1,2,3
 \eq
from which it follows that
\beq
\|V^{(k)}(\ro,\cdot)\|_{L^2(H_\ro)}\leq C\varepsilon \ro^{C\varepsilon^2},\quad k=0,1,2,3.
\eq
It remains to deduce from this a weak decay estimate. We have
\begin{lemma}\label{decay}
$$
\|V(\ro,\cdot)\|_{L^\infty(H_\ro)}^2\leq \|V(\ro,\cdot)\|_{L^2(H_\ro)} \|\pa_y V(\ro,\cdot)\|_{L^2(H_\ro)}
$$
\end{lemma}
 \begin{proof}
 By H\"older's inequality $V^2\leq 2\int |V||V_y| \, dy\leq
  2\int |V|e^{|y|/2}\, |V_y|e^{|y|/2} \, dy$.
\end{proof}
It therefore follows that
\beq
\|V^{(k)}(\ro,\cdot)\|_{L^{\infty}(H_\ro)}\leq C\ro^{C\varepsilon^2},\quad k=0,1,2.
\eq

\section{The completeness}
We now want to apply the proof of Lemma \ref{decayode} to
$$\pa_\ro^2
V+\Big(1+\frac{\beta}{\ro}V^2+\frac{1}{4\ro^2}\Big)V=F=\frac{1}{\ro^{2}}\pa_y^2
V $$ where we have proven that $|\pa_y^2 V|\leq C\varepsilon
\ro^{C\varepsilon^2}$. It then follows from the proof of that
lemma that the following limit exists $$ \Big|\big(
V_\ro^2+V^2\big)^{1/2}-a(y)\Big|\leq C\varepsilon
\ro^{-1+C\varepsilon^2}, \qquad a(y)=\lim_{\ro\to\infty} \big(
V_\ro(\ro,y)^2+V(\ro,y)^2\big)^{1/2} $$

Let
$$
 V_{\pm}=e^{\mp i\ro} \big(\pa_\ro V\pm i V\big)
$$
Then $V=(e^{i\ro} V_+-e^{-i\ro} V_-)/(2i)$, $V_+ V_-= V_\ro^2+V^2=|V_+|^2=|V_-|^2$, $V_-=\overline{V}_+$,
and
$$
 V^3=-\frac{1}{8i}\big( e^{3i\ro} V_+^3-e^{-3i\ro}V_-^3-3e^{i\ro} V_+ V_- V_+
+3e^{-i\ro} V_+ V_- V_-\big) $$
 The equation $$ \pa_\ro^2
V+V+\Big(\frac{\beta}{\ro}V^2+\frac{1}{4\ro^2}\Big)V=F $$ become
$$ \pa_\ro V_\pm+ \Big(\frac{\beta}{\ro}e^{\mp
i\ro}V^3+\frac{1}{4\ro^2}e^{\mp i\ro}V\Big) =e^{\mp i\ro }F $$ In
other words $$ \pa_\ro V_\pm\mp i g V_\pm =F_\pm+e^{\mp i\ro }F,
$$ where $$ g=\Big(\frac{\beta}{\ro}\frac{3}{8} V_+ V_-
+\frac{1}{\ro^2} \frac{1}{8} \Big) $$ $$ F_\pm=\pm
\frac{\beta}{8i\ro}\big( e^{\pm 2i\ro} V_\pm^3-e^{\mp
4i\ro}V_\mp^3 +3e^{\mp 2i\ro} V_+ V_- V_\mp \big)\mp\frac{ie^{\mp
2i\ro}}{8\ro^2} V_\mp $$ Multiplying by the integrating factor
$e^{\mp iG}$, where $G(\ro,y)=\int g(\ro,y)\, d\ro$ we get $$
\frac{d}{d\ro} \big( V_\pm e^{\mp i G}\big) =e^{\mp iG}
F_\pm+e^{\mp iG\mp i\ro} F $$ Note that $$ e^{\mp iG} F_\pm
=\frac{d}{d\ro}\Big(\frac{e^{\mp iG}}{\ro} H_\pm\Big)
+\frac{1}{\ro^2} K_\pm+\frac{1}{\ro} L_\pm F $$ where $$ H_\pm =
\pm \frac{\beta}{8i}\big( \frac{e^{\pm 2i\ro}}{\pm 2i}  V_\pm^3
-\frac{e^{\mp 4i\ro}}{\mp 4i}V_\mp^3 +3\frac{e^{\mp 2i\ro}}{\mp
2i} V_+ V_- V_\mp \big) $$ and $$ |K_\pm|\les
|V_\pm|\,(1+|V_\pm|^4),\qquad |L_\pm|\les |V_\pm|^2 $$ Hence $$
\frac{d}{d\ro}\Big( V_\pm e^{\mp iG}\mp \frac{e^{\mp iG}}{\ro}
H_{\pm} \Big)= \frac{1}{\ro^2} K_\pm +(e^{\mp iG\mp
i\ro}+\frac{1}{\ro} L_\pm) F $$ Since we have already shown that
$|V_+|=|V_-|\leq C\varepsilon$ and that $|F|\leq \ro^{-2}|\pa_y^2
V|\leq \varepsilon \ro^{-2+C\varepsilon^2}$ it follows that the
right hand side is integrable and that $$ \big| V_\pm e^{\mp
iG}-a_\pm(y)\big|\leq C\varepsilon \ro^{-1+C\varepsilon^2},
\qquad\text{where}\quad a_\pm(y)=\lim_{\ro\to\infty}
V_\pm(\ro,y)\, e^{\mp iG(\ro,y)} $$ Moreover since $V_+ V_-\sim
a(y)^2$ it follows from what we have already shown that $$ \Big|
G(\ro,y)-\beta \ln{|\ro|}\,\, \frac{3}{8}\,  a(y)^2 \Big|\leq
C\ro^{-1+C\varepsilon^2} $$ and furthermore we must have that
$|a_\pm(y)|=a(y)$. The exact phase is hence determined by looking
at $a_\pm(y)/a(y)$.

Alternatively we can use:
\begin{lemma}\label{completelemma}  Suppose that $g$ is real valued and
$$
i\pa_\ro V -g V=F
 $$
Then with $G(\ro,y)=\int_{s_0}^\ro g(\tau,y)\, d\tau$
$$
\big|\pa_\ro |V|\big| +\big|\pa_\ro (Ve^{iG})\big|\leq 2|F|\,
$$
\end{lemma}
\begin{proof} Multiplying with $\overline{V}$ gives
$$
i\pa_\ro |V|^2=\Im F \overline{V}
$$
and it follows that $|\pa_\ro |V||\leq |F|$.
Multiplying with the integrating factor $e^{iG(s)}$, where
$G=\int g\, ds$ gives
$\pa_\ro \big( Ve^{iG}\big)=-i \, F\, e^{iG}$ and the lemma follows.
\end{proof}
{\em Acknowledgement}

We would like to thank E. Taflin for discussions and comments on
an earlier version of this work.

\end{document}